\documentclass{article}
\usepackage{amsthm}
\usepackage{amsmath}
\usepackage[colorlinks=true]{hyperref}
\hypersetup{urlcolor=blue, citecolor=green}
\usepackage{todonotes}
\usepackage{amssymb}
\usepackage{xcolor}
\usepackage{tikz-cd}
\usetikzlibrary{patterns}
\usepackage{comment}
\usepackage{graphicx} % Required for inserting images
\usepackage{cite}

\newtheorem{theorem}{Theorem}
\newtheorem{proposition}{Proposition}
\newtheorem{lemma}{Lemma}

\newtheorem{remark}{Remark}

\newtheorem{coro}{Corollary}
\newtheorem{corollary}[coro]{Corollary}
\newtheorem{example}{Example}

\newtheorem{definition}{Definition}

\title{Natural Expansiveness and Equilibrium States for Semiflows}
\author{Nelda Tamblay\thanks{Departamento de Matemáticas, Facultad de
Ciencias, Universidad de Chile, Las Palmeras 3425, Ñuñoa, Santiago, Chile.
E-mail: \texttt{nelda.jaque@uchile.cl}}
\and Juan Francisco Peña\thanks{Departamento de Matemáticas, Facultad de
Ciencias, Universidad de Chile, Las Palmeras 3425, Ñuñoa, Santiago, Chile.
E-mail: \texttt{juan.pena@ug.uchile.cl}}}

\begin{document}

\maketitle

\begin{abstract}
We introduce \emph{natural expansiveness}, a notion of expansiveness for
non-invertible dynamical systems which is defined by using the natural
extension. We show that every expansive homeomorphism or flow, every
positive expansive endomorphism and every positive expansive semiflow
without end points is naturally expansive. Moreover, we prove that the
suspension of a naturally expansive endomorphism is naturally expansive.
This gives examples of naturally expansive semiflows which are neither
positive expansive nor restrictions of expansive flows. Then we extend
natural expansiveness and the topological pressure to semiflows which are
not necessarily continuous, by using a pseudosemimetric $d_\rho$. We
introduce the notion of $d_\rho$-conjugacy, and we prove that every
regular impulsive semiflow is $d_\rho$-conjugate to a continuous
semiflow. As an application, we obtain the existence and uniqueness of
equilibrium states for regular impulsive semiflows.
\end{abstract}



\section{Introduction}

The notion of expansive dynamical system was introduced by Utz
\cite{Utz} for homeomorphisms and by Bowen and Walters \cite{BW72} for
flows. Examples of expansive systems are the Anosov and Axiom A
diffeomorphisms \cite{zbMATH03321950,zbMATH03336696}. Expansiveness is a
key hypothesis in the thermodynamic formalism. For instance, Franco
\cite{Franco77} proved that an expansive flow with specification has a
unique equilibrium state for every potential with the Bowen property.
Moreover, expansiveness is invariant by conjugacy, and so these results
can be transferred between conjugate systems.

For non-invertible systems, that is, endomorphisms and semiflows, the
usual notion is positive expansiveness
\cite{zbMATH03314005,zbMATH03556927,zbMATH03762871,ACS17,Art14,zbMATH07559908}.
However, for semiflows this notion is very rigid. Herrero and Jaque
\cite{zbMATH07559908} proved that, on a compact metric space, a positive
expansive semiflow in the sense of Alves, Carvalho and Siqueira
\cite{ACS17} is the identity on a finite set. They also proved that a
positive expansive semiflow in their sense, which follows Artigue
\cite{Art14}, is a finite union of closed orbits, unbranched tails and
isolated singularities. Hence, under positive expansiveness, results for
expansive flows, such as Franco's theorem, can only be applied to very
simple semiflows.

In this article we focus on the notion of expansiveness for
non-invertible systems. We propose a notion which is not rigid, which
agrees with the classical one for invertible systems, and which can be
transferred by conjugacy. For this, we use the natural extension, that
is, the invertible system induced on the inverse limit (see for instance
\cite{zbMATH03330892,zbMATH04109690,zbMATH04078314}). We say that a
system is \emph{naturally expansive} if its natural extension is
expansive. Note that, if two systems are conjugate, then their natural
extensions are conjugate, and hence natural expansiveness is invariant
by conjugacy.

We show that natural expansiveness contains the classical notions.
Indeed, every expansive homeomorphism or flow, every positive expansive
endomorphism and every positive expansive semiflow without end points is
naturally expansive. The converse is false: the weakly Anosov endomorphisms
of Ma\~n\'e and Pugh \cite{zbMATH03493662} are naturally expansive, but in
general they are not positive expansive. Moreover, we prove that natural expansiveness recovers two
important results for flows. First, the suspension of a naturally
expansive endomorphism is naturally expansive, which extends a result of
Bowen and Walters \cite[Th.~6]{BW72}. Second, a naturally expansive
semiflow with natural specification has a unique equilibrium state, which
extends Franco's theorem \cite{Franco77}.

Finally, we consider impulsive semiflows, which are not continuous. In
\cite{BerNel19} it was proved that a regular impulsive semiflow $\phi$ is
semiconjugate to a continuous semiflow $\tilde\varphi$, by a map $H$
whose inverse is not continuous. Also, a pseudosemimetric $d_\rho$ was
introduced in order to define the topological entropy of $\phi$. In this
article we use $d_\rho$ to extend natural expansiveness, specification
and the topological pressure to semiflows which are not necessarily
continuous. Then we introduce the notion of \emph{$d_\rho$-conjugacy},
which preserves these three notions, and we prove that $H$ is a
$d_\rho$-conjugacy. As a consequence, we obtain the existence and
uniqueness of equilibrium states for regular impulsive semiflows. Note
that Alves, Carvalho and Siqueira \cite{ACS17} obtained this result
assuming positive expansiveness, which, as we saw above, is very rigid.

The article is organized as follows. In Section~2 we give the setting
and the statements of our results. In \S\ref{sec:preliminaries-exp} we
define natural expansiveness for continuous semiflows, we compare it with
the classical notions (Proposition~\ref{t1} and Theorem~\ref{t2}), and we
state the results on suspensions (Theorem~\ref{t3} and
Corollary~\ref{c1}) and on equilibrium states
(Theorem~\ref{cor:continuous-package}). In \S\ref{sec:noncont} we extend
these notions to semiflows which are not necessarily continuous, by using
$d_\rho$, and we state the results for regular impulsive semiflows
(Theorems~\ref{prop:bridge} and~\ref{thm:impulsive-equilibrium}). In
Section~\ref{sec:proofs} we give the proofs, in the same order.

\section{Setting and statements}

Throughout the article, $(X,d)$ is a compact metric space. Since we
work with natural extensions, the endomorphisms $T$ and the continuous
semiflows $\varphi$ are surjective, that is, $T(X)=X$ and
$\varphi_t(X)=X$ for all $t\ge0$.

\subsection{Expansiveness of continuous semiflows}\label{sec:preliminaries-exp}

A homeomorphism $T: X \to X$ is \textit{expansive} \cite{Utz} if there exists $\delta > 0$ such that, for all $x, y \in X$,
\[
d(T^n x, T^n y) < \delta\ \ \forall n \in \mathbb{Z}\ \Rightarrow\ x = y.
\]

Let $C_0(\mathbb{R})$ be the set of continuous maps $\gamma:\mathbb{R}\to\mathbb{R}$ with $\gamma(0)=0$. A flow $\varphi: \mathbb{R} \times X \to X$ is \textit{expansive} \cite{BW72} if for every $\varepsilon > 0$ there exists $\delta > 0$ such that, for all $x,y\in X$ and $\gamma\in C_0(\mathbb{R})$,
\[
d(\varphi_t(x), \varphi_{\gamma(t)}(y)) < \delta\ \ \forall t \in \mathbb{R}\ \Rightarrow\ \exists\, |t| < \varepsilon:\ \varphi_t(x) = y.
\]

Now we consider non-invertible systems. An endomorphism $T: X \to X$ is \textit{positive expansive} \cite{zbMATH03314005} if there exists $\delta > 0$ such that, for all $x,y\in X$,
\[
d(T^n x, T^n y) < \delta\ \ \forall n \in \mathbb{N}\ \Rightarrow\ x = y.
\]

Let $C_0(\mathbb{R}^+_0)$ be the set of continuous maps $\gamma:\mathbb{R}^+_0\to\mathbb{R}^+_0$ with $\gamma(0)=0$. A continuous semiflow $\varphi: \mathbb{R}_0^+ \times X \to X$ is \textit{positive expansive} \cite{zbMATH07559908} if for every $\varepsilon > 0$ there exists $\delta > 0$ such that, for all $x,y\in X$ and $\gamma\in C_0(\mathbb{R}^+_0)$,
\[
d(\varphi_t(x), \varphi_{\gamma(t)}(y)) < \delta\ \ \forall t \geq 0\ \Rightarrow\ \exists\, t \in [0, \varepsilon):\ \varphi_t(x) = y \text{ or } \varphi_t(y) = x.
\]

\subsubsection{Inverse limit}\label{LimiteInverso}

The \textit{inverse limit} of a system is the set of its complete orbits,
and on this set the system induces an invertible map, called the
\textit{natural extension}. Hence we can use the classical notion of
expansiveness.

Let $T:X\to X$ be an endomorphism and let $\varphi$ be a continuous
semiflow on $X$. The inverse limits of $T$ and $\varphi$ are
\[
\begin{aligned}
\overline{X}_T&=\{\overline x=(x_k)_{k\in\mathbb{Z}}:\ T(x_k)=x_{k+1}\ \text{for all } k\in\mathbb{Z}\},\\
\overline{X}_\varphi&=\{\overline{x}=(x(s))_{s\in\mathbb{R}}:\ \varphi_t(x(s))=x(s+t)\ \text{for all } t\ge0,\ s\in\mathbb{R}\},
\end{aligned}
\]
with the metrics
\[
\overline{d}_T(\overline x,\overline y)=\sum_{k\in\mathbb{Z}}\frac{d(x_k,y_k)}{2^{|k|}}
\quad\text{and}\quad
\overline{d}_\varphi(\overline{x},\overline{y})=\sum_{k\in\mathbb{Z}}\frac{d(x(k),y(k))}{2^{|k|}}.
\]
The maps $\overline{T}(\overline x)=(x_{k+1})_{k\in\mathbb Z}$ and
$\overline{\varphi}_t(\overline{x})=(x(s+t))_{s\in\mathbb{R}}$,
$t\in\mathbb R$, define a homeomorphism and a continuous flow on these
compact metric spaces (see \cite{zbMATH03330892,zbMATH02026957}). They
are called the \emph{natural extensions} of $T$ and $\varphi$. Moreover,
the projections $\pi_0(\overline x)=x_0$ and $\tau_0(\overline x)=x(0)$
are continuous semiconjugacies.

\begin{definition}[Natural expansiveness]\label{def:ile-classical}
We say that an endomorphism $T$ (resp.\ a continuous semiflow $\varphi$)
is \emph{naturally expansive} if its natural extension $\overline{T}$
(resp.\ $\overline{\varphi}$) is expansive.
\end{definition}

Note that natural expansiveness is invariant by conjugacy. Indeed, if $h$
is a conjugacy between $T$ and $S$, then
$\overline x\mapsto(h(x_k))_{k\in\mathbb Z}$ is a conjugacy between
$\overline T$ and $\overline S$. The same holds for semiflows.

\begin{remark}\label{rem:anosov}
Weakly Anosov endomorphisms are naturally expansive. Indeed, by
\cite[Prop.~2]{zbMATH03493662}, a weakly Anosov endomorphism $T$ is
uniformly expansive in the sense of Mañé and Pugh. This means that there
exists $\varepsilon>0$ such that two complete orbits of a map close to
$T$ which are $\varepsilon$-close are equal. Taking the map $T$, if
$\overline d_T(\overline T^{\,n}\overline x,\overline T^{\,n}\overline y)<\varepsilon$
for all $n\in\mathbb Z$, then $d(x_n,y_n)<\varepsilon$ for all $n$, and so
$\overline x=\overline y$.
\end{remark}

The first result is the following.

\begin{proposition}\label{t1}
Every expansive homeomorphism or flow is naturally expansive.
\end{proposition}
Following \cite{zbMATH07559908}, a point $x\in X$ is an \emph{end point}
of a semiflow $\varphi$ if $\varphi_t(x)=x$ for all $t\ge0$ and there exist
$y\neq x$ and $t>0$ such that $\varphi_t(y)=x$. Note that a flow has no end
points, since each $\varphi_t$ is injective. Moreover, the fixed points of
an expansive flow are isolated \cite{BW72}, while an end point $e$ of
$\varphi$ gives a fixed point $(e)_{s\in\mathbb R}$ of $\overline\varphi$
which is not isolated, since the complete orbits in the tail of $e$
converge to it. For this reason, for semiflows we assume that there are no end points.

\begin{theorem}\label{t2}
Every positive expansive endomorphism, and every positive expansive
semiflow without end points, is naturally expansive.
\end{theorem}

The converse of Theorem~\ref{t2} is not true. The following example shows
that there exist naturally expansive endomorphisms which are not positive
expansive.

\begin{example}
Let $T$ be the endomorphism of the torus $\mathbb R^2/\mathbb Z^2$ induced
by $A=\begin{pmatrix} n&1\\1&1\end{pmatrix}$, $n\ge3$. Then $T$ is weakly
Anosov, and so it is naturally expansive by Remark~\ref{rem:anosov}.
However, $A$ has an eigenvalue $0<\lambda_-<1$, and so the points $x$ and
$x+sv_-$, where $v_-$ is an eigenvector of $\lambda_-$ and $s$ is small,
show that $T$ is not positive expansive.
\end{example}
\subsubsection{Suspension semiflows}\label{sec:suspensions}

Let $T:X\to X$ be an endomorphism. 
The \emph{suspension} of $T$ under the constant $1$ is the semiflow
$\varphi^T$ on $X^T=X\times[0,1]/\!\sim$, where $(x,1)\sim(Tx,0)$, given by
\[
\varphi^T_t(x,u)=(T^{[u+t]}x,\{u+t\}).
\]
We consider on $X^T$ the Bowen--Walters metric $d^T$ (see
\cite[Ch.~2]{Barr13}): $d^T(z,w)$ is the infimum of
$\sum_{i=0}^{n-1}\rho^T(z_i,z_{i+1})$ over all chains
$z=z_0,z_1,\dots,z_n=w$ in which each pair $z_i,z_{i+1}$ is joined by a
horizontal or a vertical segment, where
$\rho^T((x,u),(y,u))=(1-u)d(x,y)+ud(Tx,Ty)$ for horizontal segments and
$\rho^T((x,u),(x,v))=|u-v|$ for vertical segments.

The following result shows that the inverse limit and the suspension are
compatible.

\begin{theorem}\label{t3}
Let $T$ be an endomorphism. Then $\overline{\varphi^T}$ is conjugate to
$\varphi^{\overline T}$.
\end{theorem}

Since expansiveness of flows is invariant by conjugacy and the suspension
of an expansive homeomorphism is an expansive flow \cite{BW72}, we obtain
the following corollary.

\begin{coro}\label{c1}
The suspension of a naturally expansive endomorphism is a naturally expansive semiflow.
\end{coro}
Corollary~\ref{c1} gives new examples of naturally expansive semiflows.
The following example shows that there exist naturally expansive
semiflows which are neither positive expansive nor restrictions of
flows.

\begin{example}
Let $S^1=\mathbb R/\mathbb Z$ and $T([x])=[2x]$. Since $T$ is positive
expansive, by Theorem~\ref{t2} and Corollary~\ref{c1} its suspension $\varphi^T$ is
naturally expansive. However, the points $([0],\tfrac12)$ and
$([\tfrac12],\tfrac12)$ arrive at $([0],0)$ at the time $\tfrac12$. Hence
$\varphi^T$ is not injective, and so it is not positive expansive
\cite[Prop.~3.1]{zbMATH07559908} and it is not the restriction of a flow.
\end{example}

\subsubsection{Equilibrium states for continuous semiflows}\label{sec:eq-continuous}

Let $\varphi$ be a continuous semiflow on $X$.
We denote by $M(\varphi)$ the set of $\varphi$-invariant Borel probability
measures on $X$. Given $\varepsilon,T>0$, a set $E\subset X$ is
$(\varphi,T,\varepsilon)$-separated if for all $x\ne y$ in $E$ there exists
$t\in[0,T]$ such that $d(\varphi_t(x),\varphi_t(y))\ge\varepsilon$.

Let $f:X\to\mathbb R$ be continuous. We put $S_Tf(x):=\int_0^Tf(\varphi_t(x))\,dt$
and
\[
P_T(\varphi,f,\varepsilon):=\sup\Big\{\sum_{x\in E}e^{S_Tf(x)}:
E \text{ is finite and }(\varphi,T,\varepsilon)\text{-separated}\Big\}.
\]
The \emph{topological pressure} of $f$ is defined by
$P(\varphi,f):=\lim_{\varepsilon\to0^+}P(\varphi,f,\varepsilon)$, where
$P(\varphi,f,\varepsilon):=\limsup_{T\to\infty}\frac1T\log P_T(\varphi,f,\varepsilon)$.
By the variational principle (see \cite{Ruelle,WaltersPressure}), we have
\begin{equation}\label{eq:var-principle}
P(\varphi,f)=\sup_{\mu\in M(\varphi)}\Big(h_\mu(\varphi_1)+\int_X f\,d\mu\Big).
\end{equation}

We say that a flow $\varphi$ satisfies \emph{specification} if for every
$\varepsilon>0$ there exists $L>0$ such that, for all $n\in\mathbb N$,
$x_0,\dots,x_n\in X$ and $t_0<t_1<\dots<t_{n+1}$: if $t_{i+1}-t_i\ge L$ for
all $i$, then there exist a periodic point $y\in X$, with period in
$[t_{n+1}-t_0-n\varepsilon,\,t_{n+1}-t_0+n\varepsilon]$, and a step
function $s(t)=\sum_{i}\alpha_i\chi_{[t_i,t_{i+1})}(t)$ with
$|\alpha_i-\alpha_{i+1}|<\varepsilon$, such that
\begin{equation}
d\big(\varphi_{t+s(t)}(y),\varphi_{t-t_i}(x_i)\big)<\varepsilon \quad\text{for all } i \text{ and } t\in[t_i,t_{i+1}-L].\label{def:specification}
\end{equation}
As in Definition~\ref{def:ile-classical}, we say that a
continuous semiflow $\varphi$ satisfies \emph{natural specification} if
its natural extension $\overline\varphi$ satisfies specification.

We denote by $V(\varphi)$ the set of
continuous functions $f:X\to\mathbb R$ for which there exist
$\varepsilon>0$ and $K>0$ such that, for all $t>0$ and $x,y\in X$,
\[
d(\varphi_s(x),\varphi_s(y))<\varepsilon\ \ \forall s\in[0,t]\ \Rightarrow\ \Big|\int_0^t[f(\varphi_s(x))-f(\varphi_s(y))]\,ds\Big|<K.
\]

We will use the following result of Franco \cite{Franco77}, which we call
\emph{Franco's theorem}: if $\varphi$ is an expansive flow which
satisfies specification, then for every $f\in V(\varphi)$ there exists a
unique $\mu\in M(\varphi)$ such that
$P(\varphi,f)=h_\mu(\varphi_1)+\int_X f\,d\mu$. This measure is called
the \emph{equilibrium state} of $f$. In \cite{Franco77} it is assumed
that $X$ has finite dimension, which holds by \cite{KS81}.
The following result extends
Franco's theorem to continuous semiflows by using natural expansiveness.

\begin{theorem}\label{cor:continuous-package}
Let $\varphi$ be a naturally expansive continuous semiflow which
satisfies natural specification. Then for every
$f\in V(\varphi)$ there exists a unique $\mu\in M(\varphi)$ such that
\[
P(\varphi,f)=h_\mu(\varphi_1)+\int_X f\,d\mu.
\]
\end{theorem}

\subsection{Expansiveness of semiflows not necessarily continuous}\label{sec:noncont}

A \emph{semiflow not necessarily
continuous} on $X$ is a map $\phi:\mathbb{R}_0^+\times X\to X$ such that
$\phi_0=\mathrm{id}_X$, $\phi_{t+s}=\phi_t\circ\phi_s$ for all $t,s\ge0$,
and for every $x\in X$ the orbit $t\mapsto\phi_t(x)$ has finitely many
discontinuities in each bounded interval. For instance, continuous
semiflows and impulsive semiflows (defined in \S\ref{sec:impulsive})
satisfy this condition.

The classical notions defined with the metric $d$ do not work well for
these semiflows. For instance, for an impulsive semiflow, two close points
may arrive at $D$ at different times, and so their orbits are far during
a short time. For
the topological entropy, this problem was studied in \cite{BerNel19}.
As in that work, for $\rho>0$ we consider the pseudosemimetric
$d_\rho:X\times X\to\mathbb{R}_0^+$ given by
\begin{equation}\label{def:pseudo}
d_\rho(x,y)=\inf\{d(\phi_s(x),\phi_s(y)):s\in[0,\rho)\}.
\end{equation}
The idea is that, if two orbits are far only for a short time, then
$d_\rho$ is small. Note that
\begin{equation}\label{eq:drho}
d_\rho\le d \quad\text{and}\quad d_{\rho'}\ge d_\rho \ \text{ if } \rho'<\rho,
\end{equation}
and in general $d_\rho$ does not satisfy the triangle inequality.

\subsubsection{Natural expansiveness and pressure}\label{sec:il-noncont}

The \emph{inverse limit} of $\phi$ is the set
\[
\overline{X}_\phi=\{\overline{x}=(x(s))_{s\in\mathbb{R}}:\ \phi_t(x(s))=x(s+t)\ \text{for all } t\ge0,\ s\in\mathbb{R}\},
\]
and the map $\overline\phi_t(\overline{x})=(x(s+t))_{s\in\mathbb{R}}$, $t\in\mathbb{R}$, defines a
flow (not necessarily continuous) on $\overline{X}_\phi$. Note that this definition only uses that
$\phi$ is a semiflow, and if $\phi=\varphi$ is continuous, then we obtain
the inverse limit $\overline{X}_\varphi$ of \S\ref{LimiteInverso}. For
$\rho>0$ and $\overline{x},\overline{y}\in\overline{X}_\phi$, we define
\[
\overline{d}_{\phi,\rho}(\overline{x},\overline{y}):=\sum_{k\in\mathbb{Z}}\frac{d_\rho(x(k),y(k))}{2^{|k|}}.
\]
In general, $\overline d_{\phi,\rho}$ is not a metric, but we only use it
to measure how close two orbits are.

\begin{definition}[$d_\rho$-natural expansiveness and $d_\rho$-specification]\label{def:rho-ile}
We say that $\phi$ is \emph{$d_\rho$-naturally expansive} if for every
$\varepsilon>0$ there exist $\rho>0$ and $\delta>0$ such that, for all
$\overline{x},\overline{y}\in\overline{X}_\phi$ and $\gamma\in C_0(\mathbb{R})$,
\[
\overline{d}_{\phi,\rho}\big(\overline\phi_t(\overline{x}),\overline\phi_{\gamma(t)}(\overline{y})\big)<\delta\ \ \forall t\in\mathbb{R}\ \Rightarrow\ \exists\, |t_0|<\varepsilon:\ \overline\phi_{t_0}(\overline{x})=\overline{y}.
\]

We say that $\phi$ satisfies \emph{$d_\rho$-specification} if for every
$\varepsilon>0$ there exists $\rho_0>0$ such that, for every
$\rho\in(0,\rho_0)$, there exists $L>0$ such that, for all $n\in\mathbb N$,
$\overline{x}_0,\dots,\overline{x}_n\in\overline X_\phi$ and $t_0<t_1<\dots<t_{n+1}$: if
$t_{i+1}-t_i\ge L$ for all $i$, then there exist a periodic point
$\overline{y}\in\overline X_\phi$, with period in
$[t_{n+1}-t_0-n\varepsilon,\,t_{n+1}-t_0+n\varepsilon]$, and a step
function $s(t)$ as in \eqref{def:specification}, such that
\[
\overline{d}_{\phi,\rho}\big(\overline\phi_{t+s(t)}(\overline{y}),\overline\phi_{t-t_i}(\overline{x}_i)\big)<\varepsilon \quad\text{for all } i \text{ and } t\in[t_i,t_{i+1}-L].
\]
\end{definition}

\begin{remark}\label{rem:monotone}
The name $d_\rho$ refers to the family $(d_\rho)_{\rho>0}$, since $\rho$
depends on $\varepsilon$. By \eqref{eq:drho}, if the expansiveness
condition holds for $\rho$, then it holds for every $\rho'\in(0,\rho)$
with the same $\delta$. This is not the case for specification, since
$d_\rho$ appears in the conclusion. For this reason, we ask specification
for every $\rho\in(0,\rho_0)$.
\end{remark}

The following result shows that Definition~\ref{def:rho-ile} agrees with
Definition~\ref{def:ile-classical} when the semiflow is continuous.

\begin{proposition}\label{prop:coherence}
Let $\varphi$ be a continuous semiflow. Then $\varphi$ is
$d_\rho$-naturally expansive (resp.\ satisfies $d_\rho$-specification) if
and only if $\varphi$ is naturally expansive (resp.\ satisfies natural
specification).
\end{proposition}

Now we extend the topological pressure to $\phi$ not
necessarily continuous. Let $f:X\to\mathbb R$ be a bounded Borel function.
Since the orbits of $\phi$ have finitely many discontinuities in each
bounded interval, the function $s\mapsto f(\phi_s(x))$ is Lebesgue
measurable for every $x\in X$. We put $S_Tf(x):=\int_0^Tf(\phi_s(x))\,ds$. A set
$E\subset X$ is $(\phi,T,\varepsilon,\rho)$-separated if for all $x\ne y$
in $E$ there exists $t\in[0,T]$ such that
$d_\rho(\phi_t(x),\phi_t(y))\ge\varepsilon$.

We put
\[
\hat P_T(\phi,f,\varepsilon,\rho):=\sup\Big\{\sum_{x\in E}e^{S_Tf(x)}:
E\text{ is finite and }(\phi,T,\varepsilon,\rho)\text{-separated}\Big\}
\]
and
\[
\hat P(\phi,f,\varepsilon,\rho):=\limsup_{T\to\infty}\frac1T\log\hat P_T(\phi,f,\varepsilon,\rho).
\]
Note that the function $\varepsilon\mapsto\hat P(\phi,f,\varepsilon,\rho)$
is nonincreasing. Now, we define
\[
\hat P(\phi,f,\rho):=\lim_{\varepsilon\to0^+}\hat P(\phi,f,\varepsilon,\rho).
\]
By \eqref{eq:drho}, the function $\rho\mapsto\hat P(\phi,f,\rho)$ is also
nonincreasing. We define the \emph{topological pressure with the
pseudosemimetric} by
\[
\hat P(\phi,f):=\lim_{\rho\to0^+}\hat P(\phi,f,\rho).
\]
For $f=0$, $\hat P(\phi,0)=\hat h_s(\phi)$ is the entropy defined in
\cite{BerNel19}.

The following result shows that $\hat P$ agrees with $P$ when the
semiflow is continuous.

\begin{proposition}\label{thm:pressure-coherence}
Let $\varphi$ be a continuous semiflow. Then
$\hat P(\varphi,f)=P(\varphi,f)$ for every continuous function
$f:X\to\mathbb R$.
\end{proposition}

\subsubsection{Regular impulsive semiflows}\label{sec:impulsive}

Let $\varphi:\mathbb R_0^+\times X\to X$ be a continuous semiflow on $X$,
let $D\subset X$ be a proper closed set and
let $I:D\to X$ be a continuous function. For $x\in X$ we put
\[
\tau_1(x)=
\begin{cases}
\inf\{t>0:\varphi_t(x)\in D\}, & \text{if } \varphi_t(x)\in D \text{ for some } t>0;\\
+\infty, & \text{otherwise}.
\end{cases}
\]
If $\tau_1(x)<\infty$, let $x^1=\varphi_{\tau_1(x)}(x)$. Inductively, if
$\tau_n(x)<\infty$, we put
\[
\tau_{n+1}(x)=\tau_n(x)+\tau_1(I(x^n)),\qquad x^{n+1}=\varphi_{\tau_1(I(x^n))}(I(x^n)).
\]
We say that $(X,\varphi,D,I)$ is an \emph{impulsive dynamical system} if
$\tau_1(x)>0$ and $\sup_n\tau_n(x)=\infty$ for all $x\in X$. In this case,
the \emph{impulsive semiflow} $\phi$ is defined by $\phi_t(x)=\varphi_t(x)$
for $t\in[0,\tau_1(x))$, and
\[
\phi_t(x)=\varphi_{t-\tau_n(x)}(I(x^n))\quad\text{for } t\in[\tau_n(x),\tau_{n+1}(x)).
\]
See \cite{Bon07,BBCC17,AC14,ACS17,BerNel19,TSV24} for more
details.

For $\eta>0$, we put $D_\eta=\bigcup_{x\in D}\{\varphi_t(x):0<t<\eta\}$ and
$X_\eta=X\setminus(D\cup D_\eta)$. Following \cite[Def.~2.1]{BerNel19}, we
say that $\phi$ is a \emph{regular impulsive semiflow} if
$I(D)\cap D=\emptyset$, $I$ is Lipschitz and there exist $\eta>0$ and
$0<\xi<\eta/4$ such that:
\begin{itemize}
\item[(P1)] $D_\xi$ is open;
\item[(P2)] $\varphi_\xi(D_\xi)\subset X_\xi$;
\item[(P3)] $\varphi_t(x)\notin I(D)$ for all $x\in I(D)$ and $t\in(0,\xi]$.
\end{itemize}
Note that $\phi$ is not continuous.

In \cite[Prop.~5.4 and Prop.~5.5]{BerNel19} it was proved that there exist
a compact metric space $(Y,\tilde d)$, a continuous semiflow
$\tilde\varphi$ on $Y$ and a bijection $H:X_\xi\to Y$ with
$\tilde d(H(x),H(y))\le d(x,y)$, such that
$\tilde\varphi_t\circ H=H\circ\phi_t$ for all $t\ge0$. More precisely,
$H=\pi|_{X_\xi}$ and $Y=\pi(X_\xi)$, where $\pi:X\to X/\!\sim$ is the
quotient map given by $x\sim y$ if $x=y$, $y=I(x)$, $x=I(y)$ or
$I(x)=I(y)$.

Since $H$ is a bijection, every function $f$ on $X_\xi$ induces the
function $f\circ H^{-1}$ on $Y$. However, $H^{-1}$ is not continuous, so
$f\circ H^{-1}$ is not necessarily continuous. We consider the following
class of functions.

\begin{definition}\label{def:CI}
We denote by $V_I(\phi)$ the set of bounded Borel functions
$f:X\to\mathbb R$ such that
\begin{itemize}
\item[(i)] $f$ is continuous on $X_\xi\cup D$;
\item[(ii)] $f(x)=f(I(x))$ for every $x\in D$;
\item[(iii)] there exist $\rho,\varepsilon,K>0$ such that, for all $t>0$
and $x,y\in X$,
\[
d_\rho(\phi_s(x),\phi_s(y))<\varepsilon\ \ \forall s\in[0,t]\ \Rightarrow\ \Big|\int_0^t[f(\phi_s(x))-f(\phi_s(y))]\,ds\Big|<K.
\]
\end{itemize}
\end{definition}

\begin{proposition}\label{prop:CI}
Let $\phi$ be a regular impulsive semiflow. If $f\in V_I(\phi)$, then
$f\circ H^{-1}\in V(\tilde\varphi)$.
\end{proposition}

\subsubsection{Equilibrium states for regular impulsive semiflows}\label{sec:eq-impulsive}

In order to apply Franco's theorem, we need to pass from $\phi$ to the
continuous semiflow $\tilde\varphi$. Recall that $H^{-1}$ is not
continuous. However, $H$ and $H^{-1}$ are uniformly continuous if we use
$d_\rho$ on $X_\xi$. This motivates the following definition.

\begin{definition}\label{def:drho-conj}
Let $\phi$ be a semiflow (not necessarily continuous) on $X$, let
$Z\subset X$ be a $\phi$-invariant set, and let $\psi$ be a continuous
semiflow on a compact metric space $(W,d_W)$. A bijection $h:Z\to W$ is a
\emph{$d_\rho$-conjugacy} between $\phi|_Z$ and $\psi$ if
$h\circ\phi_t=\psi_t\circ h$ for all $t\ge0$, and:
\begin{itemize}
\item[(a)] for every $\gamma>0$ there exist $\rho_1>0$ and $\delta>0$ such
that, for all $\rho\in(0,\rho_1)$ and $x,y\in Z$,
\[
d_\rho(x,y)<\delta\ \Rightarrow\ d_W(h(x),h(y))<\gamma;
\]
\item[(b)] for all $\gamma>0$ and $\rho>0$ there exists $\delta>0$ such
that, for all $x,y\in Z$,
\[
d_W(h(x),h(y))<\delta\ \Rightarrow\ d_\rho(x,y)<\gamma.
\]
\end{itemize}
\end{definition}

The following is one of the main results of this article. It shows that
every regular impulsive semiflow is $d_\rho$-conjugate to a continuous
semiflow, and that, through this $d_\rho$-conjugacy, $d_\rho$-natural
expansiveness, $d_\rho$-specification and $\hat P$ correspond to natural
expansiveness, natural specification and $P$. In particular, it
reformulates \cite[Lemma~5.6]{BerNel19}.

\begin{theorem}\label{prop:bridge}
Let $\phi$ be a regular impulsive semiflow. Then $H:X_\xi\to Y$ is a
$d_\rho$-conjugacy between $\phi|_{X_\xi}$ and $\tilde\varphi$. Moreover:
\begin{itemize}
\item[(1)] $\phi$ is $d_\rho$-naturally expansive if and only if
$\tilde\varphi$ is naturally expansive;
\item[(2)] $\phi$ satisfies $d_\rho$-specification if and only if
$\tilde\varphi$ satisfies natural specification;
\item[(3)] $\hat P(\phi,f)=P(\tilde\varphi,f\circ H^{-1})$ for every
$f\in V_I(\phi)$.
\end{itemize}
\end{theorem}

Finally, we state the main result of this section. It extends Franco's
theorem to regular impulsive semiflows. As before, $M(\phi)$ denotes the
set of $\phi$-invariant Borel probability measures on $X$.

\begin{theorem}\label{thm:impulsive-equilibrium}
Let $\phi$ be a regular impulsive semiflow. Assume that $\phi$ is $d_\rho$-naturally expansive and satisfies
$d_\rho$-specification. Then $M(\phi)\neq\emptyset$ and, for every
$f\in V_I(\phi)$,
\[
\hat P(\phi,f)=\sup_{\mu\in M(\phi)}\Big(h_\mu(\phi_1)+\int_X f\,d\mu\Big).
\]
Moreover, the supremum is attained at a unique $\mu\in M(\phi)$.
\end{theorem}

\begin{corollary}
Under the hypotheses of Theorem~\ref{thm:impulsive-equilibrium}, $\phi$
has a unique measure of maximal entropy.
\end{corollary}
Indeed, $f=0$ belongs to $V_I(\phi)$, and $\hat P(\phi,0)$ is the
entropy $\hat h_s(\phi)$ of \cite{BerNel19}.

\begin{example}\label{ex:ACS}
Let $\varphi_c$ be the suspension flow of the full shift $\sigma$ on
$\Sigma_2$ under the ceiling $c(x)=a$ if $x_0=0$ and $c(x)=b$ if $x_0=1$,
where $a,b>3$ are irrational and linearly independent over $\mathbb Q$.
Let $D=\Sigma_2\times\{1\}$ and $I(x,1)=(R(x),3)$, where $R$ reverses each
digit \cite[\S7.2]{ACS17}. Since $R$ is an isometry and $\varphi_c$ moves
the points vertically, the induced impulsive semiflow $\phi$ is regular.
Moreover, by \cite[\S7.2]{ACS17}, the restriction of $\tilde\varphi$ to
$\Lambda=\bigcap_{t\ge0}\tilde\varphi_t(Y)$ is conjugate to the suspension
flow of $\sigma$ under $c\circ R-2$, which is expansive \cite[Th.~6]{BW72}
and satisfies specification. By Proposition~\ref{t1} and
Theorem~\ref{prop:bridge}, $\phi$ is $d_\rho$-naturally expansive and
satisfies $d_\rho$-specification, and so
Theorem~\ref{thm:impulsive-equilibrium} applies to $\phi$. Note that, by
\cite[Th.~1.1]{zbMATH07559908}, positive expansiveness in the sense of
\cite{ACS17} does not hold for this example.
\end{example}

\section{Proofs}\label{sec:proofs}

In this section we prove the results of Section~2, in the same order.

\subsection{Continuous semiflows}

\subsubsection{Proof of Proposition \ref{t1} and Theorem \ref{t2}}

\begin{proof}[Proof of Proposition~\ref{t1}]
Since expansiveness is invariant by conjugacy, we only need to prove
that a homeomorphism $T$ is conjugate to $\overline T$, and that a flow
$\varphi$ is conjugate to $\overline\varphi$.

Let $T$ be a homeomorphism. For every $x\in X$, the sequence
$(T^k(x))_{k\in\mathbb Z}$ is the unique element of $\overline X_T$ with
$x_0=x$. Hence $\pi_0$ is a continuous bijection between compact metric
spaces, and so it is a homeomorphism. Since
$\pi_0\circ\overline T=T\circ\pi_0$, we conclude that $\pi_0$ is a
conjugacy between $\overline T$ and $T$.

The case of a flow $\varphi$ is analogous. We replace
$(T^k(x))_{k\in\mathbb Z}$ by $(\varphi_s(x))_{s\in\mathbb R}$ and $\pi_0$
by $\tau_0$.
\end{proof}

\begin{proof}[Proof of Theorem~\ref{t2}]
Let $T$ be a positive expansive endomorphism with constant $\delta>0$, and
let $\overline x,\overline y\in\overline X_T$ be such that
$\overline d_T(\overline T^{\,j}(\overline x),\overline T^{\,j}(\overline y))<\delta$
for all $j\in\mathbb Z$. Taking the term $k=0$ of this sum, we have
$d(x_j,y_j)<\delta$ for all $j\in\mathbb Z$. Then, for every
$k\in\mathbb Z$ and $n\in\mathbb N$,
\[
d(T^n(x_k),T^n(y_k))=d(x_{k+n},y_{k+n})<\delta,
\]
and so $x_k=y_k$. Hence $\overline x=\overline y$, that is, $\overline T$
is expansive.

We now turn to the semiflow case. Let $\varphi$ be a positive expansive
semiflow without end points. By \cite[Th.~1.2]{zbMATH07559908}, $X$ is a
union of finitely many closed orbits, unbranched tails and isolated
singularities. Since a tail is the pre-orbit of an end point, $X$ is a
union of finitely many closed orbits and isolated singularities, and so
$\varphi_t$ is injective for every $t\ge0$.

Given $\varepsilon>0$, let $\delta>0$ be the constant of positive
expansiveness of $\varphi$. By \cite{zbMATH04212997}, it is enough to
consider increasing and surjective $\gamma\in C_0(\mathbb R)$. Let
$\overline x,\overline y\in\overline X_\varphi$ be such that
$\overline d_\varphi(\overline\varphi_t(\overline x),\overline\varphi_{\gamma(t)}(\overline y))<\delta$
for all $t\in\mathbb R$. Taking the term $k=0$ and $t\ge0$, we have
\[
d(\varphi_t(x(0)),\varphi_{\gamma(t)}(y(0)))<\delta\quad\text{for all } t\ge0.
\]
Since $\gamma|_{\mathbb R_0^+}\in C_0(\mathbb R^+_0)$, there exists
$t_0\in[0,\varepsilon)$ such that $\varphi_{t_0}(x(0))=y(0)$ or
$\varphi_{t_0}(y(0))=x(0)$. Assume the first case. If $s\ge0$, then
\[
\varphi_{t_0}(x(s))=\varphi_s(\varphi_{t_0}(x(0)))=\varphi_s(y(0))=y(s).
\]
If $s<0$, then
\[
\varphi_{-s}\big(\varphi_{t_0}(x(s))\big)=\varphi_{t_0}(x(0))=y(0)=\varphi_{-s}(y(s)),
\]
and since $\varphi_{-s}$ is injective, $\varphi_{t_0}(x(s))=y(s)$. Thus $\overline y=\overline\varphi_{t_0}(\overline x)$. Analogously, in
the second case we obtain $\overline x=\overline\varphi_{t_0}(\overline y)$.
In both cases $\overline y=\overline\varphi_t(\overline x)$ for some
$|t|<\varepsilon$, which proves that $\overline\varphi$ is expansive.
\end{proof}

\subsubsection{Proof of Theorem \ref{t3}}\label{sec:t3-proof}

\begin{proof}[Proof of Theorem~\ref{t3}]
Since $\overline T$ is a homeomorphism, we write the elements of
$\overline X_T^{\overline T}$ as pairs $(\overline x,r)$, with
$\overline x\in\overline X_T$ and $r\in[0,1)$. Our candidate for the
conjugacy is the map
$h:\overline X_T^{\overline T}\to\overline{(X^T)}_{\varphi^T}$ given by
$h(\overline x,r)=(z(s))_{s\in\mathbb R}$, where
\[
z(s)=\big(x_{[s+r]},\{s+r\}\big),
\]
and $[\cdot]$ and $\{\cdot\}$ denote the integer and fractional parts.
Note that $z(k)=(x_k,r)$ for all $k\in\mathbb Z$, since $r\in[0,1)$.

First we prove that $h$ is a bijection, that is, for every
$(z(s))_{s\in\mathbb R}\in\overline{(X^T)}_{\varphi^T}$ there exists a
unique $(\overline x,r)\in\overline X_T^{\overline T}$ such that
$h(\overline x,r)=(z(s))_{s\in\mathbb R}$. Indeed, let
$z(s)=(x(s),\alpha(s))$. Since $\varphi^T$ translates the second
coordinate and $(x,1)\sim(T(x),0)$, we have $\alpha(k)=\alpha(0)$ and
$T(x(k))=x(k+1)$ for all $k\in\mathbb Z$. Hence, putting $r=\alpha(0)$
and $x_k=x(k)$, we have $\overline x=(x_k)_{k\in\mathbb Z}\in\overline X_T$
and
$h(\overline x,r)=(z(s))_{s\in\mathbb R}$. Moreover, if
$h(\overline y,q)=(z(s))_{s\in\mathbb R}$, then $(y_k,q)=z(k)=(x_k,r)$
for all $k\in\mathbb Z$, and so $(\overline y,q)=(\overline x,r)$. Note that both flows translate the second coordinate, and so
\[
h\circ\varphi^{\overline T}_t=\overline{\varphi^T}_t\circ h\quad\text{for all } t\in\mathbb R.
\]

Now, since both spaces are compact, in order to prove that $h$ is a
homeomorphism it is enough to prove that $h$ is continuous. Let $\varepsilon>0$ and $\delta=\varepsilon/6$. Let
$(\overline x,r),(\overline y,q)\in\overline X_T^{\overline T}$ be such
that $\overline d_T^{\,\overline T}((\overline x,r),(\overline y,q))<\delta$.
Then there exists a chain $z_i=(\overline z_i,p_i)$, $i=0,\dots,n$, from
$(\overline x,r)$ to $(\overline y,q)$, with
$\overline z_i=(z_{i,k})_{k\in\mathbb Z}$, such that
\[
\sum_{i=0}^{n-1}\rho^{\overline T}(z_i,z_{i+1})<2\delta.
\]
For each $k\in\mathbb Z$, the points $(z_{i,k},p_i)$ form a chain in
$X^T$ from $(x_k,r)$ to $(y_k,q)$. Hence
\[
\begin{aligned}
\overline{(d^T)}_{\varphi^T}(h(\overline x,r),h(\overline y,q))&=\sum_{k\in\mathbb Z}\frac{d^T((x_k,r),(y_k,q))}{2^{|k|}}\\
&\le\sum_{i=0}^{n-1}\sum_{k\in\mathbb Z}\frac{\rho^T((z_{i,k},p_i),(z_{i+1,k},p_{i+1}))}{2^{|k|}}.
\end{aligned}
\]
Now we fix $i$, and we put $S_i=\sum_{k\in\mathbb Z}\rho^T((z_{i,k},p_i),(z_{i+1,k},p_{i+1}))/2^{|k|}$. We have two cases.

\emph{Case 1. $z_i$ and $z_{i+1}$ form a vertical segment.} Then
\[
S_i=\sum_{k\in\mathbb Z}\frac{|p_{i+1}-p_i|}{2^{|k|}}=3\rho^{\overline T}(z_i,z_{i+1}).
\]

\emph{Case 2. $z_i$ and $z_{i+1}$ form a horizontal segment.} Since
$T(z_{i,k})=z_{i,k+1}$, we have
\[
S_i=(1-p_i)\sum_{k\in\mathbb Z}\frac{d(z_{i,k},z_{i+1,k})}{2^{|k|}}+p_i\sum_{k\in\mathbb Z}\frac{d(z_{i,k+1},z_{i+1,k+1})}{2^{|k|}}=\rho^{\overline T}(z_i,z_{i+1}).
\]

In both cases $S_i\le3\rho^{\overline T}(z_i,z_{i+1})$. Therefore
\[
\overline{(d^T)}_{\varphi^T}(h(\overline x,r),h(\overline y,q))\le3\sum_{i=0}^{n-1}\rho^{\overline T}(z_i,z_{i+1})<6\delta=\varepsilon,
\]
which proves that $h$ is continuous.
\end{proof}

\subsubsection{Proof of Theorem \ref{cor:continuous-package}}

\begin{proof}[Proof of Theorem~\ref{cor:continuous-package}]
Recall that $\tau_0\circ\overline{\varphi}_t=\varphi_t\circ\tau_0$ for all $t\ge0$.
Moreover, $d(\tau_0(\overline{x}),\tau_0(\overline{y}))\le\overline{d}_\varphi(\overline{x},\overline{y})$, since the
left side is the term $k=0$ of the sum. Hence, if $f\in V(\varphi)$, then
$f\circ\tau_0\in V(\overline{\varphi})$ with the same constants. Since
$\overline\varphi$ is expansive and satisfies specification, by Franco's
theorem there exists a unique
$\overline\nu\in M(\overline{\varphi})$ such that
\[
P(\overline{\varphi},f\circ\tau_0)=h_{\overline\nu}(\overline{\varphi}_1)+\int f\circ\tau_0\,d\overline\nu.
\]
By \cite[Th.~2.1, Prop.~3.1]{zbMATH02026957}, the map
$(\tau_0)_*:M(\overline{\varphi})\to M(\varphi)$ is a bijection which
preserves the entropy. Moreover,
$\int f\circ\tau_0\,d\overline\nu=\int f\,d(\tau_0)_*\overline\nu$. By the
variational principle \eqref{eq:var-principle}, we have
$P(\overline{\varphi},f\circ\tau_0)=P(\varphi,f)$. Hence
$\mu=(\tau_0)_*\overline\nu$ is the unique measure in $M(\varphi)$ such that
$P(\varphi,f)=h_\mu(\varphi_1)+\int f\,d\mu$.
\end{proof}

\subsection{Semiflows not necessarily continuous}

\subsubsection{Proof of Propositions \ref{prop:coherence} and \ref{thm:pressure-coherence}}

The proofs of Propositions~\ref{prop:coherence}
and~\ref{thm:pressure-coherence} follow the same argument as
\cite[Lemmas~3.1 and~3.2]{BerNel19}.

\begin{proof}[Proof of Proposition~\ref{prop:coherence}]
First, assume that $\varphi$ is $d_\rho$-naturally expansive. Let
$\varepsilon>0$ and let $\rho,\delta>0$ be as in
Definition~\ref{def:rho-ile}. Let $\overline{x},\overline{y}\in\overline X_\varphi$ and
$\gamma\in C_0(\mathbb R)$ such that
\[
\overline{d}_\varphi(\overline\varphi_t(\overline{x}),\overline\varphi_{\gamma(t)}(\overline{y}))<\delta
\quad\text{for all } t\in\mathbb R.
\]
By \eqref{eq:drho}, $d_\rho\le d$, and so
$\overline{d}_{\varphi,\rho}\le\overline{d}_\varphi$. Hence
\[
\overline{d}_{\varphi,\rho}(\overline\varphi_t(\overline{x}),\overline\varphi_{\gamma(t)}(\overline{y}))<\delta
\quad\text{for all } t\in\mathbb R.
\]
Since $\varphi$ is $d_\rho$-naturally expansive, we have
$\overline\varphi_{t_0}(\overline{x})=\overline{y}$ for some $|t_0|<\varepsilon$. Therefore
$\overline\varphi$ is expansive, that is, $\varphi$ is naturally expansive.

Now, we prove the converse. For this, first we find an upper bound of
$\overline d_\varphi$ in terms of $\overline d_{\varphi,\rho}$. Since
$\varphi$ is continuous and $X$ is a
compact metric space, for all $\alpha>0$ there exists $\beta=\beta(\alpha)>0$
such that for all $z\in X$ and $t\ge0$ we have
\[
u\in[t,t+\beta]\ \Rightarrow\ d(\varphi_t(z),\varphi_u(z))<\alpha/4.
\]
This implies that, if $\rho<\beta$, then for all $z,w\in X$ we have
\begin{equation}\label{eq:d-drho}
d(z,w)\le d_\rho(z,w)+\alpha/2.
\end{equation}
Indeed, by the triangle inequality, for all $u\in[0,\rho)$ we have
\[
d(z,w)\le d(z,\varphi_u(z))+d(\varphi_u(z),\varphi_u(w))+d(\varphi_u(w),w)
<d(\varphi_u(z),\varphi_u(w))+\alpha/2.
\]
Hence $d(z,w)-\alpha/2$ is a lower bound of
$\{d(\varphi_u(z),\varphi_u(w)):u\in[0,\rho)\}$. Since $d_\rho(z,w)$ is the
infimum of this set, we obtain \eqref{eq:d-drho}.

Then we apply \eqref{eq:d-drho} to each term of the distance in the inverse
limit. Let $\overline{x},\overline{y}\in\overline X_\varphi$. By \eqref{eq:d-drho} with
$z=x(k)$ and $w=y(k)$, and since $\sum_{k\in\mathbb Z}2^{-|k|}=3$, we have
\begin{equation}\label{eq:unif}
\overline{d}_\varphi(\overline x,\overline y)=\sum_{k\in\mathbb Z}\frac{d(x(k),y(k))}{2^{|k|}}
\le\sum_{k\in\mathbb Z}\frac{d_\rho(x(k),y(k))+\alpha/2}{2^{|k|}}
=\overline{d}_{\varphi,\rho}(\overline x,\overline y)+\tfrac{3\alpha}{2}.
\end{equation}

Assume that $\varphi$ is naturally expansive. We will prove that
$\varphi$ is $d_\rho$-naturally expansive. Let $\varepsilon>0$ and let
$\delta>0$ be the constant of expansiveness of $\overline\varphi$. Let
$\alpha=\delta/3$ and $\beta=\beta(\alpha)$, and choose $\rho\in(0,\beta)$.
Let $\overline{x},\overline{y}\in\overline X_\varphi$ and $\gamma\in C_0(\mathbb R)$ such
that
\[
\overline{d}_{\varphi,\rho}(\overline\varphi_t(\overline{x}),\overline\varphi_{\gamma(t)}(\overline{y}))<\delta/2
\quad\text{for all } t\in\mathbb R.
\]
By \eqref{eq:unif}, we have
\[
\overline{d}_\varphi(\overline\varphi_t(\overline{x}),\overline\varphi_{\gamma(t)}(\overline{y}))<\delta/2+\delta/2=\delta
\quad\text{for all } t\in\mathbb R.
\]
Since $\overline\varphi$ is expansive, we have $\overline\varphi_{t_0}(\overline{x})=\overline{y}$
for some $|t_0|<\varepsilon$. This proves that $\varphi$ is
$d_\rho$-naturally expansive, with the constants $\rho$ and $\delta/2$.

The proof for specification is analogous. If $\overline\varphi$ satisfies
specification, then, since $\overline{d}_{\varphi,\rho}\le\overline{d}_\varphi$,
$\varphi$ satisfies $d_\rho$-specification for every $\rho>0$. Conversely,
assume that $\varphi$ satisfies $d_\rho$-specification. Let
$\varepsilon>0$, let $\alpha=\varepsilon/3$ and $\beta=\beta(\alpha)$, and
let $\rho_0$ be as in Definition~\ref{def:rho-ile} for $\varepsilon/2$. We
choose $\rho<\min\{\rho_0,\beta\}$. Then the inequality
\eqref{def:specification} with $\overline{d}_{\varphi,\rho}$ and
$\varepsilon/2$ implies, by \eqref{eq:unif}, the same inequality with
$\overline{d}_\varphi$ and $\varepsilon$.
\end{proof}

\begin{proof}[Proof of Proposition~\ref{thm:pressure-coherence}]
Note that we do not use that $f\in V(\varphi)$.

First, we prove that $\hat P(\varphi,f)\le P(\varphi,f)$. By
\eqref{eq:drho}, $d_\rho\le d$. Then every
$(\varphi,T,\varepsilon,\rho)$-separated set is
$(\varphi,T,\varepsilon)$-separated. Hence
\[
\hat P_T(\varphi,f,\varepsilon,\rho)\le P_T(\varphi,f,\varepsilon),
\]
and so $\hat P(\varphi,f,\rho)\le P(\varphi,f)$ for every $\rho>0$.

Now, we prove that $P(\varphi,f)\le\hat P(\varphi,f)$. Fix $\alpha>0$ and
let $\beta=\beta(\alpha)$ be as in the proof of
Proposition~\ref{prop:coherence}. Let $E$ be a
$(\varphi,T,\alpha)$-separated set. Taking $\rho<\beta$ and
$\varepsilon\le\alpha/2$, we have that $E$ is
$(\varphi,T,\varepsilon,\rho)$-separated. Indeed, if
$d(\varphi_t(x),\varphi_t(y))\ge\alpha$, then for every $s\in[0,\rho)$ we
have
\[
d(\varphi_{t+s}(x),\varphi_{t+s}(y))\ge\alpha-\tfrac\alpha4-\tfrac\alpha4=\tfrac\alpha2\ge\varepsilon.
\]
Hence $d_\rho(\varphi_t(x),\varphi_t(y))\ge\varepsilon$.

This implies that
$P_T(\varphi,f,\alpha)\le\hat P_T(\varphi,f,\varepsilon,\rho)$. Taking
limits, we obtain
\[
P(\varphi,f,\alpha)\le\hat P(\varphi,f,\rho)\le\hat P(\varphi,f)\quad\text{for every }\rho<\beta.
\]
Taking $\alpha\to0$, we obtain the desired inequality.
\end{proof}

\subsubsection{Proof of Proposition \ref{prop:CI}}

In this subsection and the next one, $(X,\varphi,D,I)$ is a regular
impulsive dynamical system and $\phi$ is the induced impulsive semiflow.
Recall that $X_\xi\cup D$ is compact and $\pi$ is continuous
\cite[\S5]{BerNel19}.

First, we prove a kind of continuity of $H^{-1}$ when we use $d_\rho$ on
$X_\xi$. Recall that $H^{-1}$ is not continuous with the metric $d$.

\begin{lemma}\label{lem:compare}
For all $\gamma,\rho>0$ there exists $\delta>0$ such that, for all
$x,y\in X_\xi$,
\[
\tilde d(H(x),H(y))<\delta\ \Rightarrow\ d_\rho(x,y)<\gamma.
\]
\end{lemma}

\begin{proof}
We put $\tau^*(x)=\tau_1(x)$ if $x\notin D$, and $\tau^*(x)=0$ if
$x\in D$. Suppose that the lemma is not true. Then there exist
$x_n,y_n\in X_\xi$ such
that
\[
\tilde d(H(x_n),H(y_n))\to0 \quad\text{and}\quad d_\rho(x_n,y_n)\ge\gamma.
\]
Since $X_\xi\cup D$ is compact, we can assume that $x_n\to x$ and
$y_n\to y$ in $X_\xi\cup D$. Since $H=\pi$ on $X_\xi$, $\pi$ is continuous
and $\tilde d$ is a metric \cite[Prop.~5.4]{BerNel19}, we have
\[
\tilde d(\pi(x),\pi(y))=\lim_{n\to\infty}\tilde d(H(x_n),H(y_n))=0,
\]
and so $\pi(x)=\pi(y)$, that is, $x\sim y$. Since two different points
of $X_\xi$ are not equivalent \cite[\S5]{BerNel19} and
$I(D)\cap D=\emptyset$, we have three cases.

\emph{Case 1. $x=y$.}

In this case $d(x_n,y_n)\to0$. Since $d_\rho\le d$, we
have $d_\rho(x_n,y_n)\to0$, which is a contradiction.

\emph{Case 2. $x\in D$ and $y=I(x)$ (or the symmetric case).}

Since $y\in X_\xi$, we have $\tau^*(y)>0$. We fix $s\in(0,\rho)$ with
$s<\tau^*(y)$. By definition of $d_\rho$,
\[
d_\rho(x_n,y_n)\le d(\phi_s(x_n),\phi_s(y_n)).
\]
For $n$ large we have $\tau^*(x_n)<s$. Since $\tau^*(x_n)\to\tau^*(x)=0$,
$\tau^*$ is continuous on $X_\xi\cup D$ \cite[Lem.~5.2]{BerNel19}, $I$ is
continuous and $I(D)\cap D=\emptyset$, we have
$I(\varphi_{\tau^*(x_n)}(x_n))\to I(x)=y$. Hence
\[
\phi_s(x_n)=\varphi_{s-\tau^*(x_n)}\big(I(\varphi_{\tau^*(x_n)}(x_n))\big)\to\varphi_s(y).
\]
For $n$ large we have $s<\tau^*(y_n)$, and so
\[
\phi_s(y_n)=\varphi_s(y_n)\to\varphi_s(y).
\]
Since $\phi_s(x_n)$ and $\phi_s(y_n)$ converge to the same point
$\varphi_s(y)$, we have
\[
d(\phi_s(x_n),\phi_s(y_n))\to0,
\]
and so $d_\rho(x_n,y_n)\to0$, which is a contradiction.

\emph{Case 3. $x,y\in D$, $x\neq y$ and $I(x)=I(y)$.}

Let $w=I(x)=I(y)\in X_\xi$. Then $\tau^*(w)>0$. We fix $s\in(0,\rho)$ with
$s<\tau^*(w)$. As in Case 2, since $x\in D$, we have
\[
\phi_s(x_n)=\varphi_{s-\tau^*(x_n)}\big(I(\varphi_{\tau^*(x_n)}(x_n))\big)\to\varphi_s(I(x))=\phi_s(w)
\]
and, since $y\in D$,
\[
\phi_s(y_n)=\varphi_{s-\tau^*(y_n)}\big(I(\varphi_{\tau^*(y_n)}(y_n))\big)\to\varphi_s(I(y))=\phi_s(w).
\]
Since $\phi_s(x_n)$ and $\phi_s(y_n)$ converge to the same point
$\phi_s(w)$, we have
\[
d_\rho(x_n,y_n)\le d(\phi_s(x_n),\phi_s(y_n))\to0,
\]
which is a contradiction. This completes the proof.
\end{proof}

\begin{proof}[Proof of Proposition~\ref{prop:CI}]
Let $f\in V_I(\phi)$.

First, we prove that $f\circ H^{-1}$ is continuous. Let $y_n\to y$ in $Y$.
It is enough to prove that every subsequence of $(f\circ H^{-1}(y_n))$ has
a subsequence which converges to $f\circ H^{-1}(y)$. For simplicity, we
keep the notation $(y_n)$ for the subsequence.

Note that $H^{-1}(y_n)\in X_\xi$. Since $X_\xi\cup D$ is compact, there
exists $x\in X_\xi\cup D$ such that, up to a subsequence,
\[
d(H^{-1}(y_n),x)\to0.
\]
Since $\pi$ is continuous, $\pi(x)=y$. Since $f$ is continuous on
$X_\xi\cup D$ by (i), we have
\[
|f\circ H^{-1}(y_n)-f(x)|\to0.
\]
We have two cases.

\emph{Case 1. $x\in X_\xi$.}

Then $H^{-1}(y)=x$, and so $f(x)=f\circ H^{-1}(y)$.

\emph{Case 2. $x\in D$.}

Then $\pi(I(x))=\pi(x)=y$ and $I(x)\in X_\xi$, and so $H^{-1}(y)=I(x)$.
By (ii), $f(x)=f(I(x))=f\circ H^{-1}(y)$.

In both cases $f(x)=f\circ H^{-1}(y)$, and so
$f\circ H^{-1}(y_n)\to f\circ H^{-1}(y)$. This proves that
$f\circ H^{-1}$ is continuous.

Now we prove that $f\circ H^{-1}\in V(\tilde\varphi)$. Let
$\rho,\varepsilon,K$ be as in (iii), and let $\delta$ be as in
Lemma~\ref{lem:compare} with $\gamma=\varepsilon$. Let
$\tilde x,\tilde y\in Y$, and let $x=H^{-1}(\tilde x)$ and
$y=H^{-1}(\tilde y)$. Let $t>0$ and assume that
\[
\tilde d\big(\tilde\varphi_s(H(x)),\tilde\varphi_s(H(y))\big)
=\tilde d\big(H(\phi_s(x)),H(\phi_s(y))\big)<\delta
\quad\text{for all } s\in[0,t].
\]
By Lemma~\ref{lem:compare}, we have
\[
d_\rho(\phi_s(x),\phi_s(y))<\varepsilon \quad\text{for all } s\in[0,t].
\]
Hence, by (iii),
\[
\Big|\int_0^t\big[f\circ H^{-1}(\tilde\varphi_s(\tilde x))-f\circ H^{-1}(\tilde\varphi_s(\tilde y))\big]\,ds\Big|
=\Big|\int_0^t\big[f(\phi_s(x))-f(\phi_s(y))\big]\,ds\Big|<K.
\]
This proves that $f\circ H^{-1}\in V(\tilde\varphi)$.
\end{proof}

\subsubsection{Proof of Theorems \ref{prop:bridge} and \ref{thm:impulsive-equilibrium}}

First, we prove the $d_\rho$-continuity of $H$. Recall that the
$d_\rho$-continuity of $H^{-1}$ is Lemma~\ref{lem:compare}.

\begin{lemma}\label{lem:Hcont}
For every $\gamma>0$ there exist $\rho_1>0$ and $\delta>0$ such that, for
all $\rho\in(0,\rho_1)$ and $x,y\in X_\xi$,
\[
d_\rho(x,y)<\delta\ \Rightarrow\ \tilde d(H(x),H(y))<\gamma.
\]
\end{lemma}

\begin{proof}
Since $H$ and $\tilde\varphi$ are uniformly continuous (recall that
$\tilde d(H(x),H(y))\le d(x,y)$), for every $\gamma>0$ there exist
$\rho_1>0$ and $\delta>0$ such that, for all $x,y\in X_\xi$ and $\tilde z\in Y$,
\[
d(x,y)<\delta\ \Rightarrow\ \tilde d(H(x),H(y))<\gamma/3
\]
and
\[
u\in[0,\rho_1]\ \Rightarrow\ \tilde d(\tilde z,\tilde\varphi_u(\tilde z))<\gamma/3.
\]
Let $\rho\in(0,\rho_1)$ and suppose that $d_\rho(x,y)<\delta$. Then there
exists $s\in[0,\rho)$ such that
\[
d(\phi_s(x),\phi_s(y))<\delta.
\]
Hence, since $\tilde\varphi_s\circ H=H\circ\phi_s$,
\[
\begin{aligned}
\tilde d(H(x),H(y))&\le\tilde d(H(x),\tilde\varphi_s(H(x)))+\tilde d(\tilde\varphi_s(H(x)),\tilde\varphi_s(H(y)))\\
&\quad+\tilde d(\tilde\varphi_s(H(y)),H(y))\\
&<\gamma/3+\gamma/3+\gamma/3=\gamma.
\end{aligned}
\]
\end{proof}

In order to prove that $\overline H$ is a $d_\rho$-conjugacy, we use that
$\overline{(X_\xi)}_\phi=\overline X_\phi$. Indeed, it is
clear that $\overline{(X_\xi)}_\phi\subset\overline X_\phi$. Conversely,
let $\overline x\in\overline X_\phi$. By \cite[Lem.~5.2]{BerNel19}, $X_\xi$
is $\phi$-invariant and every point of $D\cup D_\xi$ arrives at $X_\xi$
before the time $2\xi$, and so $\phi_t(X)\subset X_\xi$ for all
$t\ge2\xi$. Hence, for all $s\in\mathbb R$,
\[
x(s)=\phi_{2\xi}(x(s-2\xi))\in X_\xi,
\]
and so $\overline x\in\overline{(X_\xi)}_\phi$. Therefore the map
\[
\overline H:\overline X_\phi\to\overline Y_{\tilde\varphi},\qquad \overline H(\overline{x})=(H(x(s)))_{s\in\mathbb R},
\]
is well defined.

Moreover, in order to prove the $d_\rho$-continuity of $\overline H$, note
that, since $X$ is compact, for all
$z,w\in X_\xi$ we have
\[
d_\rho(z,w)\le d(z,w)\le M\quad\text{and}\quad\tilde d(H(z),H(w))\le d(z,w)\le M,
\]
where $M$ is the diameter of $X$.

\begin{lemma}\label{lem:Hbarcont}
For every $\alpha>0$ there exist $\rho_1>0$ and $\delta>0$ such that, for
all $\rho\in(0,\rho_1)$ and $\overline{x},\overline{y}\in\overline X_\phi$,
\[
\overline d_{\phi,\rho}(\overline{x},\overline{y})<\delta\ \Rightarrow\ \overline d_{\tilde\varphi}(\overline H(\overline{x}),\overline H(\overline{y}))<\alpha.
\]
\end{lemma}

\begin{proof}
Let $\alpha>0$, and let $N$ such that $\sum_{|k|>N}M/2^{|k|}<\alpha/2$.
By Lemma~\ref{lem:Hcont} with $\gamma=\alpha/6$, there exist $\rho_1>0$ and
$\delta_1>0$ such that, for all $\rho\in(0,\rho_1)$ and $x,y\in X_\xi$,
\[
d_\rho(x,y)<\delta_1\ \Rightarrow\ \tilde d(H(x),H(y))<\alpha/6.
\]
We put $\delta=\delta_1/2^N$. Let $\rho\in(0,\rho_1)$ and suppose that
$\overline d_{\phi,\rho}(\overline x,\overline y)<\delta$. For $|k|\le N$ we
have
\[
d_\rho(x(k),y(k))\le2^{|k|}\,\overline d_{\phi,\rho}(\overline x,\overline y)<2^N\delta=\delta_1,
\]
and so $\tilde d(H(x(k)),H(y(k)))<\alpha/6$. Therefore, since $\sum_{k}2^{-|k|}=3$,
\[
\begin{aligned}
\overline d_{\tilde\varphi}(\overline H(\overline x),\overline H(\overline y))
&=\sum_{|k|\le N}\frac{\tilde d(H(x(k)),H(y(k)))}{2^{|k|}}+\sum_{|k|>N}\frac{\tilde d(H(x(k)),H(y(k)))}{2^{|k|}}\\
&<3\cdot\frac{\alpha}{6}+\frac{\alpha}{2}=\alpha.
\end{aligned}
\]
\end{proof}

Now we prove the $d_\rho$-continuity of $\overline H^{-1}$.

\begin{lemma}\label{lem:Hbarinvcont}
For all $\alpha>0$ and $\rho>0$ there exists $\delta>0$ such that, for all
$\overline{x},\overline{y}\in\overline X_\phi$,
\[
\overline d_{\tilde\varphi}(\overline H(\overline{x}),\overline H(\overline{y}))<\delta\ \Rightarrow\ \overline d_{\phi,\rho}(\overline{x},\overline{y})<\alpha.
\]
\end{lemma}

\begin{proof}
Let $\alpha>0$ and $\rho>0$, and let $N$ such that
$\sum_{|k|>N}M/2^{|k|}<\alpha/2$. By Lemma~\ref{lem:compare}, there exists $\delta_1>0$ such that, for all
$x,y\in X_\xi$,
\[
\tilde d(H(x),H(y))<\delta_1\ \Rightarrow\ d_\rho(x,y)<\alpha/6.
\]
We put $\delta=\delta_1/2^N$, and suppose that
$\overline{x},\overline{y}\in\overline X_\phi$ satisfy
$\overline d_{\tilde\varphi}(\overline H(\overline{x}),\overline H(\overline{y}))<\delta$. For
$|k|\le N$ we have
\[
\tilde d(H(x(k)),H(y(k)))\le2^{|k|}\,\overline d_{\tilde\varphi}(\overline H(\overline{x}),\overline H(\overline{y}))<2^N\delta=\delta_1,
\]
and so $d_\rho(x(k),y(k))<\alpha/6$. Hence
\[
\begin{aligned}
\overline d_{\phi,\rho}(\overline{x},\overline{y})
&=\sum_{|k|\le N}\frac{d_\rho(x(k),y(k))}{2^{|k|}}+\sum_{|k|>N}\frac{d_\rho(x(k),y(k))}{2^{|k|}}\\
&<3\cdot\frac{\alpha}{6}+\frac{\alpha}{2}=\alpha.
\end{aligned}
\]
\end{proof}

Now we prove that $\overline H$ is a $d_\rho$-conjugacy between the
natural extensions.

\begin{lemma}\label{lem:Hbarconj}
The map $\overline H:\overline X_\phi\to\overline Y_{\tilde\varphi}$ is a
$d_\rho$-conjugacy between $\overline\phi$ and $\overline{\tilde\varphi}$,
in the sense of Definition~\ref{def:drho-conj} with
$\overline d_{\phi,\rho}$ and $\overline d_{\tilde\varphi}$ instead of
$d_\rho$ and $d_W$.
\end{lemma}

\begin{proof}
By Lemmas~\ref{lem:Hbarcont} and~\ref{lem:Hbarinvcont}, it only remains to
prove that $\overline H$ is a bijection which conjugates the flows. Since
$H$ is a bijection and $H\circ\phi_t=\tilde\varphi_t\circ H$ for all
$t\ge0$, the map $\overline H$ is a bijection with inverse
$\overline H^{-1}(\overline{\tilde x})=(H^{-1}(\tilde x(s)))_{s\in\mathbb R}$.
Moreover, for all $\overline x\in\overline X_\phi$ and $t\in\mathbb R$,
\[
\overline H(\overline\phi_t(\overline x))=(H(x(s+t)))_{s\in\mathbb R}=\overline{\tilde\varphi}_t(\overline H(\overline x)).
\]
This completes the proof.
\end{proof}

Finally, using that $\overline H$ is a $d_\rho$-conjugacy
(Lemma~\ref{lem:Hbarconj}), we prove the properties which are invariant by
$d_\rho$-conjugacy.

\begin{lemma}\label{lem:inv-exp}
$\phi$ is $d_\rho$-naturally expansive if and only if $\tilde\varphi$ is
naturally expansive.
\end{lemma}

\begin{proof}
First, assume that $\phi$ is $d_\rho$-naturally expansive. Let
$\varepsilon>0$, and let $\rho,\alpha>0$ be the constants of
Definition~\ref{def:rho-ile} for $\phi$. Let $\delta$ be as in Lemma~\ref{lem:Hbarinvcont}. Let
$\overline{\tilde x},\overline{\tilde y}\in\overline Y_{\tilde\varphi}$ and $\gamma\in C_0(\mathbb R)$ such
that
\[
\overline d_{\tilde\varphi}(\overline{\tilde\varphi}_t(\overline{\tilde x}),\overline{\tilde\varphi}_{\gamma(t)}(\overline{\tilde y}))<\delta
\quad\text{for all } t\in\mathbb R.
\]
Let $\overline{x}=\overline H^{-1}(\overline{\tilde x})$ and $\overline{y}=\overline H^{-1}(\overline{\tilde y})$. Since
$\overline H$ conjugates the flows, by Lemma~\ref{lem:Hbarinvcont} we have
\[
\overline d_{\phi,\rho}(\overline\phi_t(\overline{x}),\overline\phi_{\gamma(t)}(\overline{y}))<\alpha
\quad\text{for all } t\in\mathbb R.
\]
Hence $\overline\phi_{t_0}(\overline{x})=\overline{y}$ for some $|t_0|<\varepsilon$, and so
$\overline{\tilde\varphi}_{t_0}(\overline{\tilde x})=\overline{\tilde y}$. Therefore $\tilde\varphi$ is
naturally expansive.

Now, assume that $\tilde\varphi$ is naturally expansive. Let
$\varepsilon>0$ and let $\alpha>0$ be the constant of expansiveness of
$\overline{\tilde\varphi}$. Let $\rho_1,\delta$ be as in Lemma~\ref{lem:Hbarcont}, and choose
$\rho\in(0,\rho_1)$. Let $\overline{x},\overline{y}\in\overline X_\phi$ and
$\gamma\in C_0(\mathbb R)$ such that
\[
\overline d_{\phi,\rho}(\overline\phi_t(\overline{x}),\overline\phi_{\gamma(t)}(\overline{y}))<\delta
\quad\text{for all } t\in\mathbb R.
\]
By Lemma~\ref{lem:Hbarcont}, we have
$\overline d_{\tilde\varphi}(\overline{\tilde\varphi}_t(\overline H(\overline{x})),\overline{\tilde\varphi}_{\gamma(t)}(\overline H(\overline{y})))<\alpha$
for all $t\in\mathbb R$. Hence
$\overline{\tilde\varphi}_{t_0}(\overline H(\overline{x}))=\overline H(\overline{y})$ for some
$|t_0|<\varepsilon$, and so $\overline\phi_{t_0}(\overline{x})=\overline{y}$. Therefore $\phi$
is $d_\rho$-naturally expansive.
\end{proof}

\begin{lemma}\label{lem:inv-spec}
$\phi$ satisfies $d_\rho$-specification if and only if $\tilde\varphi$
satisfies natural specification.
\end{lemma}

\begin{proof}
First, assume that $\phi$ satisfies $d_\rho$-specification. Let
$\varepsilon>0$, and let $\rho_1,\delta$ be as in Lemma~\ref{lem:Hbarcont} for
$\alpha=\varepsilon$. Let $\varepsilon'=\min\{\varepsilon,\delta\}$ and
let $\rho_0$ be as in Definition~\ref{def:rho-ile} for $\varepsilon'$. We
fix $\rho<\min\{\rho_0,\rho_1\}$, and let $L$ be the constant of the
$d_\rho$-specification for $\varepsilon'$ and $\rho$. Let
$\overline{\tilde x}_0,\dots,\overline{\tilde x}_n\in\overline Y_{\tilde\varphi}$ and $t_0<\dots<t_{n+1}$ with
$t_{i+1}-t_i\ge L$. Let $\overline{x}_i=\overline H^{-1}(\overline{\tilde x}_i)$. Then there exist a
periodic point $\overline{z}\in\overline X_\phi$, with period in
$[t_{n+1}-t_0-n\varepsilon',t_{n+1}-t_0+n\varepsilon']$, and a step
function $s(t)$ as in \eqref{def:specification} such that
\[
\overline d_{\phi,\rho}(\overline\phi_{t+s(t)}(\overline{z}),\overline\phi_{t-t_i}(\overline{x}_i))<\varepsilon'\le\delta
\quad\text{for all } i \text{ and } t\in[t_i,t_{i+1}-L].
\]
By Lemma~\ref{lem:Hbarcont} and since $\overline H$ conjugates the flows, we have
\[
\overline d_{\tilde\varphi}\big(\overline{\tilde\varphi}_{t+s(t)}(\overline H(\overline{z})),\overline{\tilde\varphi}_{t-t_i}(\overline{\tilde x}_i)\big)<\varepsilon
\quad\text{for all } i \text{ and } t\in[t_i,t_{i+1}-L].
\]
Moreover, $\overline H(\overline{z})$ is periodic with the same period as $\overline{z}$, and
$\varepsilon'\le\varepsilon$. Therefore $\overline{\tilde\varphi}$
satisfies specification.

Now, assume that $\overline{\tilde\varphi}$ satisfies specification. Let
$\varepsilon>0$ and $\rho>0$. Let $\delta$ be as in Lemma~\ref{lem:Hbarinvcont} for
$\alpha=\varepsilon$, let $\varepsilon'=\min\{\varepsilon,\delta\}$, and
let $L$ be the constant of the specification of
$\overline{\tilde\varphi}$ for $\varepsilon'$. Let
$\overline{x}_0,\dots,\overline{x}_n\in\overline X_\phi$ and $t_0<\dots<t_{n+1}$ with
$t_{i+1}-t_i\ge L$. Then there exist a periodic point
$\overline{\tilde z}\in\overline Y_{\tilde\varphi}$, with period in
$[t_{n+1}-t_0-n\varepsilon',t_{n+1}-t_0+n\varepsilon']$, and a step
function $s(t)$ as in \eqref{def:specification} such that
\[
\overline d_{\tilde\varphi}\big(\overline{\tilde\varphi}_{t+s(t)}(\overline{\tilde z}),\overline{\tilde\varphi}_{t-t_i}(\overline H(\overline{x}_i))\big)<\varepsilon'\le\delta
\quad\text{for all } i \text{ and } t\in[t_i,t_{i+1}-L].
\]
By Lemma~\ref{lem:Hbarinvcont} and since $\overline H$ conjugates the flows, the periodic point
$\overline{z}=\overline H^{-1}(\overline{\tilde z})$ satisfies
\[
\overline d_{\phi,\rho}(\overline\phi_{t+s(t)}(\overline{z}),\overline\phi_{t-t_i}(\overline{x}_i))<\varepsilon
\quad\text{for all } i \text{ and } t\in[t_i,t_{i+1}-L].
\]
Since $\rho>0$ is arbitrary, $\phi$ satisfies $d_\rho$-specification.
\end{proof}

\begin{lemma}\label{lem:inv-pressure}
For every $f\in V_I(\phi)$ we have $\hat P(\phi,f)=P(\tilde\varphi,f\circ H^{-1})$.
\end{lemma}

\begin{proof}
Let $f\in V_I(\phi)$ and $g=f\circ H^{-1}$. By Proposition~\ref{prop:CI},
$g$ is continuous. Since $f=g\circ H$ on $X_\xi$ and $X_\xi$ is
$\phi$-invariant, we have
\[
S_Tf(x)=S_Tg(H(x)) \quad\text{for every } x\in X_\xi \text{ and } T>0.
\]

First, we prove that $P(\tilde\varphi,g)\le\hat P(\phi|_{X_\xi},f)$. Let
$\varepsilon>0$, and let $\rho_1,\delta$ be as in Lemma~\ref{lem:Hcont} for
$\gamma=\varepsilon$. Let $\rho\in(0,\rho_1)$ and let $F\subset Y$ be a
$(\tilde\varphi,T,\varepsilon)$-separated set. We claim that $H^{-1}(F)$ is
$(\phi,T,\delta,\rho)$-separated. Indeed, if $x\ne y$ in $H^{-1}(F)$ and
$d_\rho(\phi_t(x),\phi_t(y))<\delta$ for all $t\in[0,T]$, then, by Lemma~\ref{lem:Hcont},
\[
\tilde d(\tilde\varphi_t(H(x)),\tilde\varphi_t(H(y)))=\tilde d(H(\phi_t(x)),H(\phi_t(y)))<\varepsilon
\quad\text{for all } t\in[0,T],
\]
which is a contradiction. Since the weights are equal, we obtain
$P_T(\tilde\varphi,g,\varepsilon)\le\hat P_T(\phi|_{X_\xi},f,\delta,\rho)$.
Taking limits, $P(\tilde\varphi,g,\varepsilon)\le\hat P(\phi|_{X_\xi},f)$,
and taking $\varepsilon\to0$ we obtain the desired inequality.

Now, we prove that $\hat P(\phi|_{X_\xi},f)\le P(\tilde\varphi,g)$. Let
$\varepsilon,\rho>0$, and let $\delta$ be as in Lemma~\ref{lem:compare} for
$\gamma=\varepsilon$. In the same way, if $E\subset X_\xi$ is
$(\phi,T,\varepsilon,\rho)$-separated, then $H(E)$ is
$(\tilde\varphi,T,\delta)$-separated. Hence
$\hat P_T(\phi|_{X_\xi},f,\varepsilon,\rho)\le P_T(\tilde\varphi,g,\delta)$.
Taking limits, we obtain the desired inequality.

Finally, $\hat P(\phi,f)=\hat P(\phi|_{X_\xi},f)$. Indeed, as in
\cite[Lemma~5.7]{BerNel19}, every point of $D\cup D_\xi$ arrives at
$X_\xi$ before the time $2\xi$. Hence, in a
$(\phi,T,\varepsilon,\rho)$-separated set, we can replace each point
outside $X_\xi$ by its first point in $X_\xi$. In this way the weights
$e^{S_Tf}$ change at most by the factor $e^{2\xi\|f\|_\infty}$, which does
not depend on $T$. Taking limits, we obtain the equality.
\end{proof}

\begin{proof}[Proof of Theorem~\ref{prop:bridge}]
Since $H$ is a bijection and $H\circ\phi_t=\tilde\varphi_t\circ H$ for all
$t\ge0$, by Lemmas~\ref{lem:Hcont} and~\ref{lem:compare}, $H$ is a
$d_\rho$-conjugacy between $\phi|_{X_\xi}$ and $\tilde\varphi$. Moreover,
(1), (2) and (3) follow from Lemmas~\ref{lem:inv-exp}, \ref{lem:inv-spec}
and~\ref{lem:inv-pressure}, respectively.
\end{proof}

\begin{proof}[Proof of Theorem~\ref{thm:impulsive-equilibrium}]
Let $f\in V_I(\phi)$ and $g=f\circ H^{-1}$. Note that the inverse limit
and the invariant measures of $\tilde\varphi$ coincide with those of
$\tilde\varphi|_\Lambda$, where $\Lambda=\bigcap_{t\ge0}\tilde\varphi_t(Y)$.
Hence we can assume that $\tilde\varphi$ is surjective.

First, we apply Theorem~\ref{cor:continuous-package} to $\tilde\varphi$.
By Theorem~\ref{prop:bridge}, $\tilde\varphi$ is naturally expansive and
satisfies natural specification. By
Proposition~\ref{prop:CI}, $g\in V(\tilde\varphi)$. Then, by
Theorem~\ref{cor:continuous-package}, there exists a unique
$\tilde\nu\in M(\tilde\varphi)$ such that
\[
P(\tilde\varphi,g)=h_{\tilde\nu}(\tilde\varphi_1)+\int g\,d\tilde\nu.
\]
By the variational principle \eqref{eq:var-principle}, we have
\[
P(\tilde\varphi,g)=\sup_{\nu\in M(\tilde\varphi)}\Big\{h_\nu(\tilde\varphi_1)+\int g\,d\nu\Big\}.
\]

Now, we pass the measures from $Y$ to $X$. Let $i:X_\xi\to X$ be the
inclusion map, and let
\[
\Psi:=(i\circ H^{-1})_*:M(\tilde\varphi)\to M(\phi).
\]
We claim that $\Psi$ is a bijection. Indeed, recall that $\phi$ is
measurable \cite{AC14} and $H$ is bimeasurable \cite{Pur66}. Since $H$ is
a bijection which conjugates $\phi|_{X_\xi}$ and $\tilde\varphi$, the map
$(H^{-1})_*$ is a bijection from $M(\tilde\varphi)$ onto the set of
$\phi$-invariant measures on $X_\xi$. Moreover, $\mu(X_\xi)=1$ for every
$\mu\in M(\phi)$ \cite[Lem.~5.1]{TSV24}, and so $i_*$ is a bijection onto
$M(\phi)$. This proves the claim. In particular, $M(\phi)\neq\emptyset$.
Since
$H\circ\phi_1=\tilde\varphi_1\circ H$ on $X_\xi$, $\Psi$ preserves the
entropy. Moreover, since $f=g\circ H$ on $X_\xi$, $\Psi$ preserves the
integral.

Finally, by Theorem~\ref{prop:bridge}(3), we have
\[
\hat P(\phi,f)=P(\tilde\varphi,g)=\sup_{\mu\in M(\phi)}\Big\{h_\mu(\phi_1)+\int f\,d\mu\Big\}.
\]
Since $\Psi$ is a bijection, $\mu:=\Psi\tilde\nu$ is the unique measure
which attains the supremum. This completes the proof.
\end{proof}

\section*{Acknowledgements}

\begin{sloppypar}
J. F. Peña was supported by ANID, Beca de Doctorado Nacional
(Folio 212110870).
\end{sloppypar}

\section*{Declaration of generative AI and AI-assisted technologies in the writing process}

During the preparation of this work, the first author used Claude
(Anthropic) to improve the language and presentation of the manuscript and
to check some details of the proofs. After using this tool, the authors
reviewed and edited the content as needed and take full responsibility for
the content of the publication.

\providecommand{\bysame}{\leavevmode\hbox to3em{\hrulefill}\thinspace}
\providecommand{\MR}{\relax\ifhmode\unskip\space\fi MR }
\providecommand{\MRhref}[2]{%
  \href{http://www.ams.org/mathscinet-getitem?mr=#1}{#2}
}
\providecommand{\href}[2]{#2}


\begin{thebibliography}{10}

\bibitem{AC14}
Jos{\'e}~F. Alves and Maria Carvalho, \emph{Invariant probability measures and
  non-wandering sets for impulsive semiflows}, Journal of Statistical Physics
  \textbf{157} (2014), no.~6, 1097--1113.

\bibitem{ACS17}
Jos\'{e}~F. Alves, Maria Carvalho, and Jaqueline Siqueira, \emph{Equilibrium
  states for impulsive semiflows}, J. Math. Anal. Appl. \textbf{451} (2017),
  no.~2, 839--857. \MR{3624769}

\bibitem{Art14}
Alfonso Artigue, \emph{Positive expansive flows}, Topology Appl. \textbf{165}
  (2014), 121--132. \MR{3174073}

\bibitem{Barr13}
Lu{\'i}s Barreira, \emph{Dimension theory of hyperbolic flows}, Springer
  Monographs in Mathematics, Springer, Cham, 2013.

\bibitem{Bon07}
Everaldo~M. Bonotto, \emph{Flows of characteristic $0^+$ in impulsive
  semidynamical systems}, J. Math. Anal. Appl. \textbf{332} (2007), no.~1,
  81--96.

\bibitem{BBCC17}
Everaldo~M. Bonotto, Matheus~C. Bortolan, Tom{\'a}s Caraballo, and Rodolfo
  Collegari, \emph{Attractors for impulsive non-autonomous dynamical systems
  and their relations}, J. Differential Equations \textbf{262} (2017), no.~6,
  3524--3550.

\bibitem{zbMATH03336696}
R.~Bowen, \emph{Periodic points and measures for axiom {A} diffeomorphisms},
  Trans. Am. Math. Soc. \textbf{154} (1971), 377--397 (English).

\bibitem{BW72}
Rufus Bowen and Peter Walters, \emph{Expansive one-parameter flows}, J.
  Differential Equations \textbf{12} (1972), 180--193. \MR{341451}

\bibitem{zbMATH03330892}
J.~R. Brown, \emph{Inverse limits, entropy and weak isomorphism for discrete
  dynamical systems}, Trans. Am. Math. Soc. \textbf{164} (1972), 55--66
  (English).

\bibitem{zbMATH03314005}
Murray Eisenberg, \emph{Expansive transformation semigroups of endomorphisms},
  Fundam. Math. \textbf{59} (1966), 313--321 (English).

\bibitem{Franco77}
Ernesto Franco, \emph{Flows with unique equilibrium states}, Amer. J. Math.
  \textbf{99} (1977), no.~3, 486--514.

\bibitem{zbMATH02026957}
Lianfa He, Lina Zhou, and Jinfeng Lu, \emph{The dynamical properties of a
  continuous semiflow and its inverse limit}, Chaos Solitons Fractals
  \textbf{14} (2002), no.~8, 1169--1177.

\bibitem{zbMATH07559908}
Sebasti{\'a}n Herrero and Nelda Jaque, \emph{On two notions of expansiveness
  for continuous semiflows}, J. Math. Anal. Appl. \textbf{515} (2022), no.~1,
  10 (English), Id/No 126405.

\bibitem{BerNel19}
Nelda Jaque and Bernardo San~Mart{\'i}n, \emph{Topological entropy for
  discontinuous semiflows}, J. Differential Equations \textbf{266} (2019),
  no.~6, 3580--3600.

\bibitem{KS81}
Harvey~B. Keynes and Mark Sears, \emph{Real-expansive flows and topological
  dimension}, Ergodic Theory and Dynamical Systems \textbf{1} (1981), no.~2,
  179--195.

\bibitem{zbMATH03556927}
Ricardo Ma\~{n}e, \emph{Axiom {A} for endomorphisms}, Geom. {Topol}., {III}.
  {Lat}. {Am}. {Sch}. {Math}., {Proc}., {Rio} de {Janeiro} 1976, {Lect}.
  {Notes} {Math}. 597, 379-388 (1977)., 1977.

\bibitem{zbMATH03493662}
Ricardo Mane and Charles Pugh, \emph{Stability of endomorphisms}, Dyn. {Syst}.,
  {Proc}. {Symp}. {Univ}. {Warwick} 1973/74, {Lect}. {Notes} {Math}. 468,
  175-184 (1975)., 1975.

\bibitem{zbMATH04212997}
Masatoshi Oka, \emph{Expansiveness of real flows}, Tsukuba J. Math. \textbf{14}
  (1990), no.~1, 1--8 (English).

\bibitem{Pur66}
Roger Purves, \emph{Bimeasurable functions}, Fundamenta Mathematicae
  \textbf{58} (1966), 149--157.

\bibitem{zbMATH04078314}
J.~Quandt, \emph{On structural stability for semiflows}, Commun. Math. Phys.
  \textbf{117} (1988), no.~2, 191--202 (English).

\bibitem{zbMATH04109690}
\bysame, \emph{On inverse limit stability for maps}, J. Differ. Equations
  \textbf{79} (1989), no.~2, 316--339 (English).

\bibitem{zbMATH03762871}
William~L. Reddy, \emph{Expanding maps on compact metric spaces}, Topology
  Appl. \textbf{13} (1982), 327--334 (English).

\bibitem{Ruelle}
David Ruelle, \emph{Statistical mechanics on a compact set with zp action
  satisfying expansiveness and specification}, Transactions of the American
  Mathematical Society \textbf{185} (1973), 237--251.

\bibitem{zbMATH03321950}
S.~Smale, \emph{Differentiable dynamical systems. {With} an appendix to the
  first part of the paper: ``{Anosov} diffeomorphisms'' by {John} {Mather}},
  Bull. Am. Math. Soc. \textbf{73} (1967), 747--817 (English).

\bibitem{TSV24}
Nelda Tamblay, Bernardo San~Mart{\'i}n, and Kendry Vivas, \emph{Topological
  entropy and metric entropy for regular impulsive semiflows}, Journal of
  Dynamics and Differential Equations \textbf{36} (2024), 3553--3568.

\bibitem{Utz}
W.~R. Utz, \emph{Unstable homeomorphisms}, Proceedings of the American
  Mathematical Society \textbf{1} (1950), no.~6, 769--774.

\bibitem{WaltersPressure}
Peter Walters, \emph{A variational principle for the pressure of continuous
  transformations}, American Journal of Mathematics \textbf{97} (1975), no.~4,
  937--971.

\end{thebibliography}
\end{document}